\documentclass[12pt]{amsart}
\usepackage{amsmath,amsfonts,amssymb,amsthm,amstext,pgf,graphicx,verbatim,lmodern,textcomp,color,young,tikz}
\usetikzlibrary{decorations}
\usepackage[mathscr]{euscript}
\usetikzlibrary{decorations.markings}
\usetikzlibrary{arrows}
\usepackage{float}
\usepackage{enumerate}
\usepackage{hyperref}
\theoremstyle{plain}
\newtheorem{theorem}{Theorem}[section]
\newtheorem{lemma}[theorem]{Lemma}
\newtheorem{proposition}[theorem]{Proposition}

\theoremstyle{remark}
\newtheorem{remark}[theorem]{Remark}

\theoremstyle{definition}
\newtheorem{defn}{Definition}[section]

\newtheorem{question}{Question}[section]
\newtheorem{problem}{Problem}[section]

\title{}
\begin{document}
\title [Difference of intersection power graph and power graph]{On the difference of the intersection power graph and the power graph of a finite group}
\author[Sudip Bera]{Sudip Bera}
\author[Peter J. Cameron]{Peter J. Cameron}
\address[Sudip Bera]{Faculty of Mathematics, Dhirubhai Ambani University, Gandhinagar, India}
\email{sudip\_bera@dau.ac.in}
\address[Peter J. Cameron]{School of Mathematics and Statistics, University of St. Andrews, UK}
\email{pjc20@st-andrews.ac.uk} 

%\subjclass[2010]{ 05C25, 05C17}
		
		\begin{abstract}
The \emph{difference graph} $\mathcal{D}(G)$ of a finite group $G$ is the graph obtained by taking the (edge)
difference of the intersection power graph and the power graph of $G$, and subsequently removing all isolated vertices. In this paper, we give a number of results about the difference graph. 
We examine groups whose power graph and intersection power graph coincide. 
In addition, we make some observations on isolated vertices in difference graphs. 
We study the connectedness and perfectness of difference graph 
with respect to various properties of the underlying group $G$. 
Furthermore, we investigate the operation of twin reduction on graphs, 
a technique that yields smaller graphs which may be easier to analyze.

\medskip
%\end{abstract}
		
%\begin{keyword}
			
			\medskip  
			
			
%\end{keyword}
\end{abstract}

\maketitle

\section{Introduction}
\label{sec:introduction}

The study of graphs associated with various algebraic structures has attracted 
considerable attention over the past two decades. Such investigations are useful in 
\begin{enumerate}
    \item characterizing the resulting graphs,
    \item characterizing algebraic structures that give rise to isomorphic graphs, and
    \item understanding the interdependence between algebraic structures and their 
    corresponding graphs.
\end{enumerate}
A wide variety of graphs have been introduced to explore the properties of 
algebraic structures through graph-theoretic methods. Notable examples include 
the commuting graph~\cite{Brauer-Fowler}, the generating graph~\cite{Guralnick-Kantor_generating-graph}, 
the power graph~\cite{undpwrgraphofsemgmainsgc1,powergraph}, and the enhanced power graph~\cite{aacns,bera-dey-JGT,Bera-Dey-Mukherjee}. 
The notion of a power graph was originally introduced in semigroup theory by 
Kelarev and Quinn~\cite{Kelarev-Quinn-2000} and the intersection power graph was introduced in \cite{intersectionpwegraphb3}.

In this paper, we define the power graph and intersection power graph of a group. For 
comparison, we also define the enhanced power graph, though we do not study
this graph.

Throughout this paper, all groups $G$ under consideration are finite. Let
$\pi(G) = \{\, p \in \mathbb{N} : p \text{ is prime and } p \mid |G| \,\}$
denote the set of prime divisors of $|G|$. Similarly, for $a \in G$, let
$\pi(a) = \{\, p \in \mathbb{N} : p \text{ is prime and } p \mid o(a) \,\}$
denote the set of prime divisors of the order of $a$.

\begin{defn}\label{Def:ipg}\cite{intersectionpwegraphb3} 
Let $G$ be a group with identity element $e$.
\begin{enumerate}
\item The \emph{intersection  power graph} of $G$, denoted by $\mathcal{G}_I(G)$,  is the graph with vertex set $G$ and two distinct non-identity vertices $x$ and $y$ are adjacent if and only if there exists a \emph{non-identity} element $z\in G$
such that $z$ is a power of both $x$ and $y$; equivalently, $\langle x\rangle\cap \langle y\rangle \neq \{e\}$. By convention, we take $e$ to be adjacent to all other vertices, where $e$ is the identity element of the group $G$.   
\item
The \emph{power graph} $\mathcal{P}(G)$ of $G$ is the graph with vertex set $G$ 
with two vertices $x$ and $y$ adjacent if and only if one of them is a power of the other.   
\item
The \emph{enhanced power graph} $\mathcal{P}_e(G)$ of $G$ is the graph with 
vertex set $G$, in which two vertices $x$ and $y$ are adjacent if there exists
$u\in G$ such that both $x$ and $y$ are powers of $u$.
\end{enumerate}
\end{defn}

\begin{remark}
The power graph has received a lot of attention in the literature: see
\cite{powergraph} for a survey. Less attention has been given to the
intersection power graph: see~\cite{intersectionpwegraphb3}. It was briefly
considered in~\cite{gog} under the name \emph{dual enhanced power graph}.
 
Clearly, if $x$ and $y$ are joined in the power graph of $G$, then they are joined in both the intersection power graph and the enhanced power graph; in other
words, the power graph is a spanning subgraph of each of the other two.
\end{remark}

This motivates the study of two difference graphs; in one of these, the edges
are those of the enhanced power graph which are not edges of the power graph;
in the other, they are those of the intersection power graph which are not
edges of the power graph. The first of these graphs was studied in \cite{bcdd},
and will not be considered further here; our purpose is to study the second.

\begin{defn}\label{Def:dg}
The \emph{undeleted difference graph} $\mathbb{D}(G)$ of a finite group $G$ is the (edge) difference of the intersection power graph and the power graph of $G$, while the \emph{difference graph} $\mathcal{D}(G)$ of $G$ is the subgraph of $\mathbb{D}(G)$ obtained by deleting all isolated vertices of $\mathbb{D}(G)$.  
\end{defn}

The graph $\mathcal{D}(G)$ can be constructed from the semilattice of cyclic subgroups of $G$ as follows.
Let $G$ be a finite group and let $L(G)$ denote the semilattice of cyclic subgroups of $G$. Let $\mathcal{B}_1(G)$ be the simple undirected graph with vertex set $L(G)$, in which two vertices $S, T \in L(G)$ are adjacent if and only if $S \not\subseteq T$, $T \not\subseteq S$, and $S \cap T \neq \{e\}$, where $e$ denotes the identity element of $G$. Let $\mathcal{B}(G)$, called the \emph{base graph}, be the graph obtained from $\mathcal{B}_1(G)$ by removing all isolated vertices. The graph $\mathcal{D}(G)$ is then obtained from $\mathcal{B}(G)$ by replacing each vertex $S$ with an independent set $V(S)$ of size $|V(S)| = \varphi(|S|)$, where $\varphi$ denotes Euler's totient function, such that $x \sim y$ in $\mathcal{D}(G)$ if and only if there exist $S \neq T$ in $L(G)$ with $x \in V(S)$, $y \in V(T)$, and $S \sim T$ in $\mathcal{B}(G)$.

\begin{comment}
It is worth noting that for each edge $S \sim T$ in $\mathcal{B}(G)$, the subgraph of $\mathcal{D}(G)$ induced by $V(S) \cup V(T)$ is complete bipartite (see Lemma~3.7). Consequently, $\mathcal{D}(G)$ is empty (respectively, connected) if and only if $\mathcal{B}(G)$ is empty (respectively, connected). Thus, for the purposes of studying emptiness or connectedness, it suffices to consider $\mathcal{B}(G)$ rather than the whole of $\mathcal{D}(G)$.

For a cyclic group $G$ of order $n$, the graph $\mathcal{B}(G)$ admits a simple arithmetical description, as follows. Let $n$ be a natural number and let $S(n)$ denote the set of all divisors of $n$. Define the graph $\mathcal{B}_1(n)$ with vertex set $S(n)$, in which two vertices $x, y \in S(n)$ are adjacent if and only if $x \nmid y$, $y \nmid x$, and $\gcd(x,y) > 1$. The graph $\mathcal{B}(n)$ is obtained from $\mathcal{B}_1(n)$ by deleting all isolated vertices.

For a cyclic group $G$ of order $n$, we then have $\mathcal{B}(G) = \mathcal{B}(n)$, and $\mathcal{D}(G)$ is obtained from $\mathcal{B}(n)$ by replacing each vertex $x$ with an independent set $V(x)$ of size $|V(x)| = \varphi(x)$, such that $a \sim b$ in $\mathcal{D}(G)$ if and only if there exist $x \neq y$ in $S(n)$ with $a \in V(x)$, $b \in V(y)$, and $x \sim y$ in $\mathcal{B}(n)$.

\end{comment}
This graph is the subject of the paper. However, the first question that should
be asked after this definition is the following:

\begin{problem}
Given a group $G$, identify the isolated vertices in the undeleted difference graph $\mathbb{D}(G)$. In
particular, which groups have the property that $\mathbb{D}(G)$ has no
edges?
\end{problem}
A \emph{null graph} is an edgeless graph with non-empty vertex
set whereas an \emph{empty graph} has empty set of vertices.  
Note that a group $G$ has the property that $\mathbb{D}(G)$
has no edges if and only if the intersection power graph and the power graph of $G$
are equal. We will discuss further the question of isolated vertices in 
Section~\ref{sec:G-is-p-grp-diff-graph-empty}, and the question of when it is
the null graph in Sections~\ref{s:null} and~\ref{sec:simple-grp}.

\begin{remark}
An analogous question, that of determining groups whose power graph and
enhanced power graph are equal, has been solved~\cite{aacns}. These are
precisely the \emph{EPPO groups}, those in which every element has prime
power order. After preliminary work by Higman and Suzuki, these were
classified by Brandl~\cite{Brandl} in 1981. (A perhaps more accessible version
appears in the survey~\cite{cm}.) We hope that the groups for 
which the power graph and intersection power graph are equal also form an
interesting class.
\end{remark}

\subsection{Plan of the Paper}

The paper is organized as follows. In this section, we introduce the power graph, the intersection power graph, and the difference graph $\mathcal{D}(G)$ associated to a finite group $G$, and we discuss the questions that motivate our study.
In Section~\ref{s:null}, we study the finite groups $G$ for which $\mathcal{D}(G)$ is empty, and we record several observations concerning isolated vertices.
Section~\ref{sec:basic-props} develops the basic theory of difference graphs. We establish adjacency criteria involving generators of $G$ and show that every finite graph arises as an induced subgraph of $\mathcal{D}(G)$ for some finite group $G$.
In Section~\ref{sec:connected}, we investigate the connectivity of $\mathcal{D}(G)$. This study is motivated by a question of Aalipour et al.~\cite[Question~40]{aacns} on the connectivity of power graphs, subsequently answered by Cameron and Jafari~\cite{cj}, and by the analogous result of Bera, Dey, and Mukherjee for enhanced power graphs~\cite{Bera-Dey-Mukherjee}.
Section~\ref{sec:G-is-p-grp-diff-graph-empty} characterizes the finite $p$-groups $G$ for which $\mathcal{D}(G)$ is empty, with particular attention to the generalized quaternion groups.
In Section~\ref{sec:perfect-bipartite-eulerian-diff-graph}, we examine perfectness, bipartiteness, and the Eulerian property of $\mathcal{D}(G)$. This investigation is motivated by the fact that power graphs are always perfect~\cite[Thm.~5]{AAA}, whereas the perfectness of enhanced power graphs remains open, despite known results on their chromatic number and weak perfection~\cite{Cameron-Veronica}.
Section~\ref{sec:simple-grp} shows that the non-abelian finite simple groups $G$ for which $\mathcal{D}(G)$ is empty are precisely those listed in Proposition~\ref{Prop:null-diff-simple-grp}.
Section~\ref{sec:grp-with-unique-involution} treats groups possessing a unique involution, and Section~\ref{sec:example} illustrates the general theory through the Mathieu group $M_{11}$. Finally, we conclude in Section~\ref{sec:Conclusion-open-prob} with a discussion of open problems, including the classification of isolated vertices and the determination of perfectness, bipartiteness, and connectivity of $\mathcal{D}(G)$ when $\pi(Z(G)) \leq 1$.

\section{When is the difference graph empty?}
\label{s:null}

In this section, we examine groups whose power graph and intersection power
graph are equal, that is, those groups whose difference graph $\mathcal{D}(G)$ is empty, or equivalently, whose undeleted difference graph $\mathbb{D}(G)$ is a null graph. We also make several observations concerning the isolated vertices of the undeleted difference graph $\mathbb{D}(G)$.
We begin with a simple result showing that the undeleted difference graph shares a
property with many other graphs on groups.

\begin{proposition}
Let $H$ be a subgroup of $G$. Then the undeleted difference graph of $H$ is
the induced subgraph of the undeleted difference graph of $G$ on the vertex
set $H$.
\end{proposition}

This is true because the adjacency of vertices in $H$ depends only on elements
of $H$.

Let $\mathcal{C}$ denote the class of finite groups whose difference graph is
empty. The next theorem raises the possibility of using induction to classify
the groups in $\mathcal{C}$.

\begin{theorem}\label{t:induct}
Suppose that either $G$ cannot be generated by two elements, or $G$ has
trivial center. Then $\mathcal{D}(G)$ is empty if and only if $\mathcal{D}(H)$
is empty for every proper subgroup $H$ of $G$.
\end{theorem}

\begin{proof}
Given two elements $x$ and $y$, if the subgroup $H=\langle x,y\rangle$ is
proper, and $x$ and $y$ are joined in $\mathcal{D}(G)$ if and only if they are
joined in $\mathcal{D}(H)$. So we only have to consider the case when
$\langle x,y\rangle=G$. But in this case $\langle x\rangle\cap\langle y\rangle$
commutes with $x$ and $y$, and so lies in $Z(G)$, and so by assumption is
trivial; so $x$ and $y$ are not joined in the intersection power graph.
\end{proof}

There is a simple sufficient condition for membership in $\mathcal{C}$:

\begin{proposition}\label{p-elts}
\begin{enumerate}[(a)]
\item An element of a group $G$ whose order is prime is isolated in the undeleted
difference graph.
\item If every element of $G$ has prime order, then $G\in\mathcal{C}$.
\end{enumerate}
\end{proposition}

\begin{proof}
(a) Suppose that $a$ has prime order $p$. If $b$ is joined to $a$ in the  undeleted
difference graph, then there exists $c\in\langle a\rangle\cap\langle b\rangle$
with $c\ne e$. Thus $\langle a\rangle=\langle c\rangle$, so
$a\in\langle b\rangle$, and $a$ is joined to $b$ in the power graph, a
contradiction.

\smallskip

(b) is clear from (a).
\end{proof}

\begin{remark}
Groups satisfying condition (b) of this Proposition are EPPO groups, so can
be identified from Brandl's list. In particular, $A_5$ is the only simple
group satisfying this condition, since the other simple groups in the list
contain elements of order $4$ or $9$. However, there are other simple groups
in $\mathcal{C}$, as we will see.
\end{remark}

We conclude this section with a curious connection with the power graph. Recall
that a graph $\Gamma$ is a \emph{cograph} if it does not contain the $4$-vertex
path $P_4$ as induced subgraph.

\begin{proposition}\label{p:nulld_implies_cograph}
If $\mathcal{D}(G)$ is empty, then the power graph of $G$ is a cograph.
\end{proposition}

\begin{proof} We show the contrapositive. Suppose that $\mathcal{P}(G)$ is
not a cograph; let $(a,b,c,d)$ be an induced $P_4$. In the directed power
graph (with an arc $x\to y$ if $y$ is a power of $x$), the directions must
alternate; for, if $a\to b\to c$, say, then $a\to c$, and $a$ and $c$ would
be joined in the power graph. Suppose, without loss, that
$a\to b\gets c\to d$. Then $b$ is a power of both $a$ and $c$; and $b\ne e$,
since $b$ is not joined to $d$. Thus $a$ and $c$ are joined in the intersection
power graph but not in the power graph; so their difference is not empty.
\end{proof}

Finite groups whose power graph is a cograph are studied in \cite{cmm}, though
these authors do not give a complete classification. We could note that, if the power graph  is a cograph, then it is chordal \cite{B-F_p-pw-co}.

\begin{remark}
The converse of this theorem is false: for $G=\mathbb{Z}_4\times\mathbb{Z}_2$, the power
graph is a cograph, but the difference graph is not empty. However, it does
hold in many cases, as we will see.
\end{remark}

\section{Basic properties}
\label{sec:basic-props}
In this section, we review some known results and establish a few lemmas that will play a crucial role in the subsequent sections. Moreover, we prove that the difference graphs are universal.

\begin{lemma}[Lemma 3.10, \cite{bera-dey-JGT}]
Let $G$ be a $p$-group. Then the number of distinct subgroups of $G$ with order~$p$ is either $1$ or at least $3$.    
\end{lemma}

\begin{lemma}\cite[Theorem 5.4.10]{gorenstein} \label{Lemma:Burnside}
If $G$ is a $p$-group with a unique subgroup of order $p$ for an odd prime $p,$ then $G$ is cyclic.    
\end{lemma}

\begin{lemma}[Theorem 3.1, \cite{intersectionpwegraphb3}]\label{ipg_complete}
Let $G$ be a group. Then the intersection power graph $\mathcal{G}_I(G)$ of the group $G$ is complete if and only if either $G$ is a cyclic $p$-group or $G$ is a generalized quaternion group.
\end{lemma}

\begin{lemma}[Theorem 2.12, \cite{undpwrgraphofsemgmainsgc1}]\label{cyc_complete}
Let $G$ be a group. Then the power graph $\mathcal{P}(G)$ of the group $G$ is complete if and only if $G$ is cyclic group of order $1$ or $p^m,$ for some prime number $p$ and for some $m\in \mathbb{N}.$
\end{lemma}

\begin{lemma}[Theorem 4, \cite{cj}]\label{l:dom}
Let \(G\) be a finite group, and suppose \(x \in G\) has the property that
 for all $y\in G$, either $x$ is a power of $y$ or vice versa. Then exactly one of the following cases holds:
\begin{enumerate}
  \item \(x = e\).
  \item \(G\) is cyclic and \(x\) is a generator.
  \item \(G\) is a cyclic \(p\)-group (for some prime \(p\)), and \(x\) is an arbitrary element.
  \item \(G\) is a generalized quaternion group, and \(x\) has order \(2\).
\end{enumerate}
\end{lemma}

Let $\mathcal{G}_a$ be the set of all generators of the cyclic group generated
by the element $a$. Then we have the following:
\begin{lemma}\label{l:gens_null}
Let $G$ be a group. Then for each non-identity element $a\in G$, the subgraph induced by $\mathcal{G}_a$ is an empty subgraph of $\mathcal{D}(G)$.
\end{lemma}

\begin{proof}
The proof follows from the facts that the vertices in $\mathcal{G}_a$ form a clique in both the intersection power graph and the power graph.   
\end{proof}

\begin{lemma}\label{l:blowup}
Let $G$ be a finite group and $\mathcal{G}_a\neq \mathcal{G}_b$ for two distinct elements $a, b \in G$. If an element of $\mathcal{G}_a$ is adjacent to an element of $\mathcal{G}_b$, then each element of $\mathcal{G}_a$ is adjacent to every element of $\mathcal{G}_b$.    
\end{lemma}

\begin{proof}
Let $\mathcal{G}_a\neq \mathcal{G}_b$. Suppose that $x\in \mathcal{G}_a$ and $y\in \mathcal{G}_b$ such that $x\sim y$ in the difference graph $\mathcal{D}(G).$ Therefore, $x\sim y$ in the intersection power graph $\mathcal{G}_I(G)$ but they are not edge connected in the power graph $\mathcal{P}(G).$ Now, $x'\sim y'$ in $\mathcal{D}(G)$ for any $x'\in \mathcal{G}_a$ and $y'\in \mathcal{G}_b$ as $\langle x\rangle=\langle x'\rangle$ and $\langle y\rangle=\langle y'\rangle.$  
\end{proof}

\subsection{Difference graphs are universal}

\begin{theorem}\label{t:univ}
For every finite graph $\Gamma$, there is a finite group $G$ such that
$\Gamma$ is an induced subgraph of $\mathcal{D}(G)$. Moreover, $G$ can be
chosen to be cyclic of squarefree order.
\end{theorem}

The proof depends on two results which may be of independent interest.

Recall that a \emph{Sperner family} is a family $\mathcal{F}$ of subsets of
a finite set $X$ with the property that, for any two distinct sets
$F_1,F_2\in\mathcal{F}$, we have $F_1\not\subseteq F_2$.

\begin{proposition}
Any finite graph is an induced subgraph of the intersection graph of a
Sperner family.
\end{proposition}

\paragraph{Proof} Let $\Gamma=(V,E)$ be a finite graph. Take $X=V\cup E$,
and for each $v\in V$, put $F_v=\{v\}\cup\{e\in E:v\in e\}$. Then:
\begin{itemize}
\item If $v\ne w$, then $F_v\not\subseteq F_w$: for $v\in F_v\setminus F_w$.
\item If $v$ is joined to $w$, then $e=\{v,w\}\in F_v\cap F_w$; otherwise
$F_v\cap F_w=\emptyset$.
\end{itemize}
So $\mathcal{F}=\{F_v:v\in V\}$ is a Sperner family and the map $v\mapsto F_v$
is an isomorphism from $\Gamma$ to its intersection graph.

\begin{proposition}
The intersection graph of any Sperner family is an induced subgraph
of $\mathcal{D}(G)$, where $G$ is a cyclic group of squarefree order.
\end{proposition}

\paragraph{Proof} Let $(X,\mathcal{F})$ be a Sperner family. Choose a prime
$p_x$ for each $x\in X$, so that the chosen primes are distinct. Now let
\[G=\prod_{x\in X}\mathbb{Z}_{p_x},\]
where $g_x$ is a generator of the factor $\mathbb{Z}_{p_x}$. Then $G$ is
a cyclic group of squarefree order. For
each set $F\in\mathcal{F}$, put $g(F)=\prod_{x\in F}g_x$. Since $\mathcal{F}$
is a Sperner family, none of the elements $g(F)$ is a power of another;
and $g(F)$ and $g(F')$ are joined in the difference graph if and only if
$F\cap F'\ne\emptyset$.

\begin{remark} 
We could instead choose the groups $G$ in the theorem to
be nonabelian simple groups. For the group $\mathrm{PSL}(2,q)$ contains a
cyclic subgroup of order $(q-1)/2$. So to embed the cyclic group of order $n$
in it, we choose $q$ to be a prime congruent to $1$ modulo $2n$; by Dirichlet's
theorem, there are infinitely many such primes.
\end{remark}

%%%%%%%%%%%%%%%%%%%%%%%%%%%%%%%%%%%%%%%%%%%%%%%%%%%%%%%%%%%%%%%%%%%%%%%%%%%%%%%%%% BIPARTITE thm for nilpotent grp %%%%%%%%%%%%%%%%%%%%%%%%%%%%%%%%%%%%%%%%%%%%%%%%%%%%%%%%%%%%%%

\section{When is the difference graph connected?}
\label{sec:connected}

In this section, we examine groups whose difference graphs are connected. In the connected case, we establish bounds on the diameter of these difference graphs, and in the disconnected case, we make observations concerning their isolated vertices. Our main theorems in this regard are as follows:

\begin{theorem}\label{t:conn}
Let $G$ be a finite group with $\pi(Z(G))\ge2.$ Then the difference graph $\mathcal D(G)$ is connected if and only if $G\not\cong\mathbb Z_{pq}$ (with $p,q$ primes). Moreover, when connected,
\[
\operatorname{diam}\bigl(\mathcal D(G)\bigr)\leq 6.
\]
\end{theorem}

Let $\mathcal G(\mathbb Z_m)$ denote the set of all generators of the cyclic group $\mathbb Z_m$, and let $\mathcal I(\mathbb{D}(G))$ denote the set of all isolated vertices in the graph $\mathbb{D}(G)$.

\begin{theorem}\label{t:isol}
Under the same hypotheses, the set of isolated vertices is
\[
\mathcal I(\mathbb{D}(G)) =
\begin{cases}
G, & \text{if } G \cong \mathbb Z_{pq},\; p\ne q\text{ primes};\\
\mathcal G(\mathbb Z_m)\cup\{e\}, & \text{if } G \cong \mathbb Z_m,\; m\ne pq;\\
\{\,g\in G : o(g)\text{ is prime}\}\cup\{e\}, & \text{otherwise}.
\end{cases}
\]
\end{theorem}

\begin{comment}
As an immediate consequence of Theorems \ref{t:conn} and \ref{t:isol}, we have the following:
\begin{theorem}\label{t:twoprimes}
Let $\pi(\mathbb{Z}_m)\geq 2.$ Then the difference graph $\mathcal{D}(\mathbb{Z}_m)$ is empty if and only if $m=pq,$ where $p$ and $q$ are two distinct primes.     
\end{theorem}
\end{comment}

\subsection{Cyclic group}
For a cyclic group $G$ of order $n$, the graph $\mathcal{B}(G)$ admits a simple arithmetical description, as follows. Let $n$ be a natural number and let $S(n)$ denote the set of all divisors of $n$. Define the graph $\mathcal{B}_1(n)$ with vertex set $S(n)$, in which two vertices $x, y \in S(n)$ are adjacent if and only if $x \nmid y$, $y \nmid x$, and $\gcd(x,y) > 1$. The graph $\mathcal{B}(n)$ is obtained from $\mathcal{B}_1(n)$ by deleting all isolated vertices. For a cyclic group $G$ of order $n$, we then have $\mathcal{B}(G) = \mathcal{B}(n)$, and $\mathcal{D}(G)$ is obtained from $\mathcal{B}(n)$ by replacing each vertex $x$ with an independent set $V(x)$ of size $|V(x)| = \varphi(x)$, such that $a \sim b$ in $\mathcal{D}(G)$ if and only if there exist $x \neq y$ in $S(n)$ with $a \in V(x)$, $b \in V(y)$, and $x \sim y$ in $\mathcal{B}(n)$.

\begin{lemma}\label{lemma:diam(B(G))-g-cyclic}
Let $G$ be a cyclic group of order $n$. Then the graph $\mathcal{B}(G) = \mathcal{B}(n)$ is empty if and only if $\pi(G) \leq 1$ or $G \cong \mathbb{Z}_{pq}$ for some distinct primes $p,q$. In all other cases, $\mathcal{B}(G)$ is connected with
\[
\operatorname{diam}\big(\mathcal{B}(G)\big) \leq 3.
\]
\end{lemma}

\begin{proof}
If $\pi(G) \leq 1$ or $G \cong \mathbb{Z}_{pq},$ then clearly $\mathcal{B}_1(G)$ is an empty graph.

Conversely,
let $x, y \mid n$ with $x, y \notin \pi(n) \cup \{1, n\}$ and $x \neq y$, where $\pi(n)$ is the set of prime divisors of a natural number $n$. We consider the following cases.

\medskip
\noindent\textbf{Case I:} $n = p_1 p_2 \cdots p_k$, $k \geq 3.$ 

\smallskip
\noindent\textit{Case IA.} If $\pi(x) \not\subseteq \pi(y)$ and $\pi(y) \not\subseteq \pi(x)$, then there exist primes $p_i, p_j$ with $p_i \mid x$, $p_i \nmid y$, $p_j \mid y$, $p_j \nmid x$. Hence
$x \sim p_ip_j \sim y.$

\smallskip
\noindent\textit{Case IB.} If $\pi(x) \subseteq \pi(y)$, then since $y \neq n$ we have $\pi(y) \neq \pi(n)$; choose $p_k \in \pi(n) \setminus \pi(y)$ and $p_i \in \pi(x)$. Then $x \sim p_ip_k \sim y.$

\medskip
\noindent\textbf{Case II:} $n = p_1^{\alpha_1}p_2^{\alpha_2}\cdots p_k^{\alpha_k}$, $k \geq 2$, with some $\alpha_i > 1$.

\smallskip
\noindent\textit{Case IIA.} Analogous to Case IA.

\smallskip
\noindent\textit{Case IIB.} If $\pi(x) \subseteq \pi(y) \subsetneq \pi(n)$, the argument is analogous to Case IB.

\smallskip
\noindent\textit{Case IIC.} Suppose $\pi(x) \subsetneq \pi(y) = \pi(n)$. Since $y \neq n$, there exists $k$ with $\alpha_k > 1$ and $p_k^{\alpha_k} \nmid y$; in particular $y \nmid x$.
\begin{itemize}
\item If $x \nmid y$, then $x \sim y$ directly.
\item If $x \mid y$ and $x = p_k^{r_k}$ with $1 < r_k < \alpha_k$, then for any $i \neq k$,
\[
x = p_k^{r_k} \sim p_ip_k \sim p_k^{\alpha_k} \sim y.
\]
\item If $x \mid y$ and there exists $p_i \mid x$ with $i \neq k$, then
\[
x \sim p_ip_k \sim p_k^{\alpha_k} \sim y \quad \text{if } p_k \nmid x, \qquad\text{and}\qquad x \sim p_k^{\alpha_k} \sim y \quad \text{if } p_k \mid x.
\]
\end{itemize}

\smallskip
\noindent\textit{Case IID.} Suppose $\pi(x) = \pi(y) = \pi(n)$.
\begin{itemize}
\item If $x \nmid y$ and $y \nmid x$, then $x \sim y$ directly.
\item If $x \mid y$ and $y \nmid x$: since $y \neq n$, there exists $k$ with $\alpha_k > 1$ and $p_k^{\alpha_k} \nmid y$, whence also $p_k^{\alpha_k} \nmid x$. Then
\[
x \sim p_k^{\alpha_k} \sim y.
\]
\end{itemize}

In every case, $d(x,y) \leq 3$, which proves the bound.
\end{proof}
\begin{remark}\label{remrk-tide-bound-B(G)}
The bound stated in Lemma \ref{lemma:diam(B(G))-g-cyclic} is sharp. For $G \cong \mathbb{Z}_{p^3q}$, the graph $B(p^3q)$ is precisely the path
\[
p^2\sim  pq\sim p^3 \sim p^2q,
\]
which has diameter exactly $3.$ 
\end{remark}

\begin{theorem}\label{t:cyclic-diam}
Let $G$ be a cyclic group of order $n$. Then $\mathcal{D}(G)$ is empty if and only if $\pi(G) \leq 1$, or $G \cong \mathbb{Z}_{pq}$ for some distinct primes $p$ and $q$. In all other cases, $\mathcal{D}(G)$ is connected, and
\[
\operatorname{diam}\big(\mathcal{D}(G)\big) \leq 3.
\]
Moreover, this bound is sharp.
\end{theorem}

\begin{proof}

The proof follows from the discussion above, Lemma \ref{lemma:diam(B(G))-g-cyclic} and Remark \ref{remrk-tide-bound-B(G)}.    
\end{proof}

\subsection{Neither $G \cong \mathbb{Z}_m$ nor $G \cong \mathbb{Z}_k \times Q_{2^n}$, where $\gcd(k,2)=1$}

\smallskip
In this subsection, we assume that the group $G$ is neither cyclic nor a direct product of a cyclic group of odd order and a generalized quaternion group. We assume
that $\pi(Z(G))\ge 2.$

Let the order of the group be given by 
\[
|G| = p_1^{\alpha_1} p_2^{\alpha_2} \cdots p_r^{\alpha_r}, \quad \text{with } r \geq 2.
\]
For each $j \in [r]$, let $\langle a_{j1} \rangle, \langle a_{j2} \rangle, \dots, \langle a_{jk_j} \rangle$ be the distinct cyclic subgroups of $G$ of order $p_j$, where $k_j \geq 1$.

Assume that for some fixed $m, n \in [r]$, there exists $a_{mx}, a_{ny} \in Z(G)$ such that $\operatorname{o}(a_{mx}) = p_m$ and $\operatorname{o}(a_{ny}) = p_n$ for some $x \in [k_m], y\in [k_n]$. For each $j \in [r] \setminus \{m\}$, define the set
\[
U_j = \{ a_{mx} a_{j1}, a_{mx} a_{j2}, \dots, a_{mx} a_{jk_j} \}. 
\]
Also, we define the set
\[
V_j = \{ a_{ny} a_{j1}, a_{ny} a_{j2}, \dots, a_{ny} a_{jk_j} \},   
\]
for each $ j\in [r] \setminus \{n\}$.
Next, define
\[
\mathcal{U}_x = \bigcup_{j \in [r] \setminus \{m\}} U_j, \quad \mathcal{V}_y = \bigcup_{j \in [r] \setminus \{n\}} V_j  \quad \text{and} \quad \mathcal{V}=\mathcal{U}_x\bigcup \mathcal{V}_y.
\]

%\Let $X = \{ a_{jt} : j \in [r] \setminus \{1\},\ t \in [k_j] \}.$

\begin{lemma}\label{t:clique}
The vertices in $\mathcal{U}_x$ ( in $\mathcal{V}_y$) form a clique in the difference graph $\mathcal{D}(G).$  
\end{lemma}
\begin{proof}
 Notice that the order of each member of $\mathcal{U}_x$ is a product of two distinct primes. Moreover, $a_{js}, a_{mx}\in \langle a_{mx}a_{js}\rangle,$ for any $j\in [r]\setminus\{m\}$ and $s\in [k_j].$
 Therefore, no two vertices in $\mathcal{U}_x$ are edge connected in the power graph $\mathcal{P}(G)$, and any two vertices in $\mathcal{U}_x$ are edge connected in the intersection power graph of $G.$ Similarly, any two vertices in $\mathcal{V}_y$ are edge connected in $\mathcal{D}(G).$ This completes the proof.    
\end{proof}
\begin{lemma}\label{t:pconn2}
The vertices in $\mathcal{V}$ are path connected by a path of length at most $2$ in $\mathcal{D}(G).$    
\end{lemma}
\begin{proof}
Let $g_1, g_2\in \mathcal{V}=\mathcal{U}_x\bigcup \mathcal{V}_y.$ If both $g_1$ and $g_2$ belong either in $\mathcal{U}_x$ or in $\mathcal{V}_y,$ then by Lemma \ref{t:clique}, they are edge connected. Now, suppose $g_1\in \mathcal{U}_x$ and $g_2\in \mathcal{V}_y.$ Then $g_1 = a_{mx}a_{qs}$ and $g_2 = a_{ny}a_{uv}.$ Note that $a_{ny}a_{mx}=a_{mx}a_{ny}$ and $a_{mx}a_{ny}\in \mathcal{V}=\mathcal{U}_x\bigcap \mathcal{V}_y.$ Therefore, by Lemma \ref{t:clique}, $g_1\sim a_{mx}a_{ny}\sim g_2$ in $\mathcal{D}(G).$  
\end{proof}

\begin{lemma}\label{t:k3}
Under the same hypotheses, the difference graph $\mathcal{D}(G)$ contains a cycle of length~$3$.
\end{lemma}

\begin{proof}
We consider two cases:

\medskip

\textbf{Case 1:} Suppose at least one of $k_m, k_n$ equals $1$, and without loss of generality, let $k_m = 1$. Since $\pi(G) \ge 2$ and $G$ is neither cyclic nor isomorphic to a direct product of a cyclic group of odd order with a generalized quaternion group, there exists some prime index $j \in [r] \setminus \{m\}$ for which $k_j \ge 3$. Choose three distinct elements $a_{j1}, a_{j2}, a_{j3}$, and let $a_{mx} \in Z(G)$. Then the vertices
\[
a_{mx}a_{j1}, \quad a_{mx}a_{j2}, \quad a_{mx}a_{j3}
\]
in $\mathcal{D}(G)$ form a $3$‑cycle, since each pair is adjacent.

\medskip

\textbf{Case 2:} If neither $k_m$ nor $k_n$ equals $1$, then $k_n \ge 3$. Again, take $a_{mx} \in Z(G)$ and three distinct elements $a_{n1}, a_{n2}, a_{n3}$ in the corresponding subgroup. The vertices
\[
a_{mx}a_{n1}, \quad a_{mx}a_{n2}, \quad a_{mx}a_{n3}
\]
similarly form a $3$‑cycle in $\mathcal{D}(G)$.

Thus, in both cases, a triangle exists in the difference graph.
\end{proof}

\begin{theorem}\label{t:pconn3}  
Let $g \in V(\mathcal{D}(G)) \setminus \mathcal{V}.$ If $\pi(Z(G)) \geq 2,$ then there exists $h \in \mathcal{V}$ such that $g$ and $h$ are path connected in $\mathcal{D}(G).$ Moreover, $d(g,h) \leq 2.$  
\end{theorem}

\begin{proof}
The order of each member of $\mathcal{V}$ is a product of two distinct primes.  
We will prove that $g$ is path connected with some $h \in \mathcal{V}$ in $\mathcal{D}(G).$  
Since $\pi(Z(G)) \geq 2,$ let $a_{ny} \in Z(G)$ such that $\mathrm{o}(a_{ny}) = p_n,$ for some $y \in [k_n]$ and $p_m \neq p_n.$  

We consider the following cases.

\medskip
\noindent
\textbf{Case 1.} Suppose $p_m \nmid \mathrm{o}(g).$  
Then there exists a prime $p_i \neq p_m$ dividing $\mathrm{o}(g).$  
Let $a_{is} \in \langle g \rangle,$ for some $s \in [k_i].$  
Consider the element $a_{mx}a_{is} \in \mathcal{V}_x.$  
Since 
\[
a_{is} \in \langle a_{mx}a_{is} \rangle \cap \langle g \rangle,
\]
we have $g \sim a_{mx}a_{is}$ in $\mathcal{G}_I(G).$  
Moreover, $\mathrm{o}(g)$ is not prime (by Proposition~\ref{p-elts}) and $p_m \nmid \mathrm{o}(g),$ so $g \nsim a_{mx}a_{is}$ in the power graph $\mathcal{P}(G).$  
Hence, $g \sim a_{mx}a_{is}$ in $\mathcal{D}(G).$

\medskip
\noindent
\textbf{Case 2.} Suppose $p_m \mid \mathrm{o}(g).$  
We consider two subcases.

\smallskip
\textbf{Subcase 1.} Let $a_{mx} \in \langle g \rangle.$  
If there exists a prime $p_t \in \{p_1,\dots, p_r\}$ such that $p_t \nmid \mathrm{o}(g),$ then take $a_{mx}a_{tz},$ where $\mathrm{o}(a_{tz}) = p_t.$  
Clearly, $a_{mx}a_{tz} \in \mathcal{V}$ and 
\[
a_{mx} \in \langle a_{mx}a_{tz} \rangle \cap \langle g \rangle.
\]
Also, neither $\langle a_{mx}a_{tz} \rangle \subseteq \langle g \rangle$ nor $\langle g \rangle \subseteq \langle a_{mx}a_{tz} \rangle.$  
Thus, $a_{mx}a_{tz} \sim g$ in $\mathcal{D}(G).$  

Now suppose that $\mathrm{o}(g)$ is divisible by all primes $p_1,\dots, p_r.$  
Since $G$ is neither cyclic nor a direct product of a cyclic group of odd order and a generalized quaternion group, there exists a prime $p_j \in \{p_1,\dots,p_r\}$ such that the number of distinct cyclic subgroups of order $p_j$ is $k_j \geq 3$ (by Lemma \ref{Lemma:Burnside}). If $p_j = p_m,$ take $a_{mu}a_{iv},$ where $\mathrm{o}(a_{mu}) = p_m,$ $\mathrm{o}(a_{iv}) = p_i \neq p_m,$ with $\langle a_{mu} \rangle \neq \langle a_{mx} \rangle$ and $a_{iv} \in \langle g \rangle,$ for some $v \in [k_i].$  
Then $g \sim a_{mu}a_{iv}$ in $\mathcal{D}(G).$ If $p_j \neq p_m,$ then $G$ has a unique cyclic subgroup $\langle a_{mx} \rangle$ of order $p_m.$  
Let $a_{jw} \in G$ with $\mathrm{o}(a_{jw}) = p_j$ and $a_{jw} \notin \langle g \rangle.$  
Then $g \sim a_{mx}a_{jw}$ in $\mathcal{D}(G).$

\smallskip
\textbf{Subcase 2.} If $a_{mx} \notin \langle g \rangle.$  
Let $a_{mv} \in \langle g \rangle$ with $\mathrm{o}(a_{mv}) = p_m.$ First, suppose $\mathrm{o}(g) = p_m^d$ with $d \geq 2.$  
Now, $a_{ny} \in Z(G)$ with $\mathrm{o}(a_{ny}) = p_n.$  
Consider $a_{mv}a_{ny}.$  
Although $a_{mv}a_{ny} \notin \mathcal{V},$ we have
\[
a_{mv} \in \langle a_{mv}a_{ny} \rangle \cap \langle g \rangle,
\]
hence $g \sim a_{mv}a_{ny}$ in $\mathcal{D}(G).$  
Also, $a_{mx}a_{ny} \in \mathcal{V}$ and $a_{mv}a_{ny} \sim a_{mx}a_{ny}$ in $\mathcal{D}(G),$ since 
\[
a_{ny} \in \langle a_{mv}a_{ny} \rangle \cap \langle a_{mx}a_{ny} \rangle.
\]
Thus, $g \sim a_{mv}a_{ny} \sim a_{mx}a_{ny}$ in $\mathcal{D}(G).$ If $\mathrm{o}(g)$ has at least two distinct prime divisors, then it follows similarly that $g \sim a_{mx}a_{sw},$ where $\mathrm{o}(a_{sw}) = p_s \neq p_m$ and $a_{sw} \in \langle g \rangle.$  

\medskip
This completes the proof.  
\end{proof}

\begin{theorem}\label{t:conn2}
Let $G$ be a finite group which is neither cyclic nor a direct product of a cyclic group of odd order and a generalized quaternion group. If $\pi(Z(G))\geq 2$, then $\mathcal{D}(G)$ is connected and $\text{diam}(\mathcal{D}(G))\leq 6$. 
\end{theorem}
\begin{proof}
The proof of this theorem follows from Lemma \ref{t:pconn2} and Theorem \ref{t:pconn3}. 
\end{proof}

\subsection{$G\cong \mathbb{Z}_m\times Q_{2^n},$ where $\gcd(m, 2)=1$}

\begin{lemma}\label{t:cxq_conn}
Let \(G \cong \mathbb{Z}_m \times Q_{2^n}\), where \(\gcd(m,2)=1\) and \(n \ge 3\). Then the difference graph \(\mathcal{D}(G)\) is connected, and moreover,
$\operatorname{diam}\bigl(\mathcal{D}(G)\bigr) \le 3.$

\end{lemma}

\begin{proof}
Write \(m = p_1^{\alpha_1} \cdots p_r^{\alpha_r}\) with each \(\alpha_i \ge 1\). For each \(i\), choose an element \(a_i \in G\) of order \(p_i\), and let \(b \in G\) be the unique element of order \(2\).
Recall that the number of distinct subgroups of order \(4\) in \(Q_{2^n}\) is \(2^{n-2} + 1\). Let \(b_1,\dots,b_{2^{n-2}+1}\) be representatives such that \(\langle b_i \rangle\) are these distinct subgroups.
Now define
\[
\mathcal{A} = \{\,a_1b, \ldots, a_rb\} \;\cup\; \{\,b_1,\ldots,b_{2^{n-2}+1}\}.
\]
Each pair of vertices in \(\mathcal{A}\) is adjacent in \(\mathcal{D}(G)\), so \(\mathcal{A}\) forms a clique. Next, let \(g \in V(\mathcal{D}(G)) \setminus \mathcal{A}\). There are two cases:

\medskip

\textbf{Case 1.} If \(\gcd\bigl(\mathrm{o}(g), 2\bigr)=1\), pick any prime divisor \(p_i\) of \(\mathrm{o}(g)\). Then \(g\) is adjacent to the vertex \(a_i b \in \mathcal{A}\).

\medskip

\textbf{Case 2.} If \(\gcd\bigl(\mathrm{o}(g), 2\bigr) \neq 1\), then again \(g\) is adjacent to each \(a_i b\) in \(\mathcal{A}\).

Since \(\mathcal{A}\) forms a clique, any two vertices in \(V(\mathcal{D}(G))\) can be joined by a path of length at most \(3\), at worst going through two vertices in \(\mathcal{A}\):
\[
g \sim a_i b \sim a_j b \sim h.
\]
This establishes the connectivity of \(\mathcal{D}(G)\) and the diameter bound.
\end{proof}

\begin{proof}[Proof of Theorem \ref{t:conn}]
Proof of this theorem follows from Theorems \ref{t:cyclic-diam}, \ref{t:conn2}, and Lemma \ref{t:cxq_conn}.
\end{proof}

\begin{proof}[Proof of Theorem \ref{t:isol}]
Suppose \(G \cong \mathbb{Z}_{m}\). Then by Lemma \ref{l:dom}, each generator of \(G\) is an isolated vertex in the difference graph \(\mathcal{D}(G)\).
The complete classification of all isolated vertices now follows directly from Proposition \ref{p-elts} together with Theorem \ref{t:conn}.
\end{proof}

%%%%%%%%%%%%%%%%%%%%%%%%%%%%%%%%%%%%%%% Isolated vertex for p-group %%%%%%%%%%%%%%%%%%%%%%%%%%%%%%%%%%%%%%%%%%%%%%%%%%%%%%%%%%%%%%%%%%%%%%%
\section{When $G$ is a $p$-group}
\label{sec:G-is-p-grp-diff-graph-empty}

In this section, we investigate finite \(p\)-groups for which the power graph and the intersection power graph coincide, equivalently, those groups whose difference graph is empty. We also discuss observations concerning isolated vertices in such difference graphs.

\begin{theorem}\label{t:empty}
Let \(G\) be a finite \(p\)-group. Then \(\mathcal{D}(G)\) is an empty graph if and only if \(G\) satisfies one of the following conditions:
\begin{enumerate}
  \item \(G\) is a cyclic \(p\)-group.
  \item \(G\) is union of at least three cyclic subgroups such that they have pairwise trivial intersection.
    
\end{enumerate}
\end{theorem}

To prepare for the proof of Theorem~\ref{t:empty}, we establish several supporting lemmas below.

\begin{lemma}\label{t:gq}
Let \(Q_{2^n}\) be the generalized quaternion group of order \(2^n\), with \(n \ge 3\). Denote by \(H = \langle x \rangle\) the unique subgroup of index 2 in \(Q_{2^n}\), i.e., of order \(2^{n-1}\). Then for any element \(g \in Q_{2^n} \setminus H\) and any \(h \in H\) with \(\operatorname{o}(h) \neq 2\), the vertices \(g\) and \(h\) are adjacent in the difference graph \(\mathcal{D}(G)\).
\end{lemma}

\begin{proof}
Let \(g \in Q_{2^n} \setminus H\). Then \(\operatorname{o}(g) = 4\), and there are exactly \(2^{n-2}\) distinct subgroups of order 4 in \(Q_{2^n}\setminus H\). Denote by \(x\) the unique element of order 2 in \(Q_{2^n}\); note that \(x \in \langle g \rangle \cap \langle h \rangle\). Therefore, \(g\) and \(h\) are adjacent in the intersection power graph \(\mathcal{G}_I(G)\). Since they are not adjacent in the power graph \(\mathcal{P}(G)\), it follows by definition that \(g\) is adjacent to \(h\) in the difference graph \(\mathcal{D}(G)\).
\end{proof}

\begin{lemma}\label{t:gq_isol}
Let \(Q_{2^n}\) be the generalized quaternion group of order \(2^n\). A non-identity element
$g \in Q_{2^n} \setminus \{e\}$ is an isolated vertex in the undeleted difference graph \(\mathbb{D}(Q_{2^n})\) if and only if \(\operatorname{o}(g) = 2\).
\end{lemma}

\begin{proof}
Let \( g \) be the unique element of order \(2\) (the involution) in \(Q_{2^n}\). By Proposition~\ref{p-elts}, any non‑identity element of prime order in a finite group is an isolated vertex of the difference graph \(\mathcal{D}(G)\). Consequently, \(g\) is isolated in \(\mathcal{D}(Q_{2^n})\).

Now consider any element \( z \in Q_{2^n} \) with order \( \operatorname{o}(z) > 2 \). By Lemma~\ref{t:gq}, such an element \(z\) is adjacent to some other element in \(Q_{2^n}\) within \(\mathcal{D}(G)\). Hence \(z\) is not isolated.

Combining these observations completes the proof.
\end{proof}

\noindent
We now focus on a finite \(p\)-group \(G\) that is neither cyclic nor generalized quaternion. In this situation, there are at least three distinct cyclic subgroups of order \(p\). Denote these by
\[
H_1 = \langle a_1 \rangle, \quad \ldots, \quad H_s = \langle a_s \rangle, 
\quad s \ge 3.
\]
For each \(i \in [s]\), define
\[
K_i = \{\,g \in G : a_i \in \langle g \rangle\,\}.
\]

\begin{lemma}\label{t:totdisc}
Let \(G\) be a finite \(p\)-group that is neither cyclic nor generalized quaternion. Then the difference graph \(\mathcal{D}(G)\) is empty if and only if
\[
G = \mathbb{Z}_{p^{t_1}} \,\cup\cdots\cup\,\mathbb{Z}_{p^{t_s}},
\quad s \ge 3.
\]
\end{lemma}

\begin{proof}
Suppose 
\[
G = \mathbb{Z}_{p^{t_1}} \cup \cdots \cup \mathbb{Z}_{p^{t_s}},
\quad s \ge 3.
\]
Since each \(\mathbb{Z}_{p^{t_i}}\) is a cyclic $p$-group. So, by Lemmas \ref{ipg_complete} and \ref{cyc_complete}, $\mathcal{D}(G)$ is empty. 

Conversely, assume \(\mathcal{D}(G)\) is an empty graph. Since \(G\) is neither cyclic nor generalized quaternion, there are at least \(s \ge 3\) distinct cyclic subgroups of order \(p\), giving rise to a decomposition
\[
G = K_1 \cup \cdots \cup K_s\cup \{e\}.
\]
We show each \(K_i\cup \{e\}\) is itself a cyclic subgroup.
Let \(x \in K_i\) be an element of maximal order among those in \(K_i\). Since \(x\) is isolated in \(\mathcal{D}(G)\), it remains adjacent to every other \(y \in K_i\) in the power graph \(\mathcal{P}(G)\). By the maximality of \(\mathrm{o}(x)\), this implies \(y \in \langle x \rangle\) for all \(y \in K_i\). Thus \(K_i\cup \{e\} = \langle x \rangle\), a cyclic group.
This completes the proof.
\end{proof}

\begin{proof}[Proof of Theorem~\ref{t:empty}]
If \(G\) is a cyclic \(p\)-group, then by Lemmas \ref{ipg_complete} and \ref{cyc_complete} no two distinct vertices are adjacent in \(\mathcal{D}(G)\), so the graph is edgeless.

For the non‑cyclic case, the conclusion follows immediately from Lemmas~\ref{t:gq_isol} and \ref{t:totdisc}. These lemmas together complete the characterization of when \(\mathcal{D}(G)\) is empty.
\end{proof}

\section{Perfectness, bipartiteness, and Eulerian properties of difference Graphs}
\label{sec:perfect-bipartite-eulerian-diff-graph}
In this section, we address the following fundamental questions for a finite group \(G\) with \(\pi(Z(G)) \ge 2\):

\begin{itemize}
  \item When is the difference graph \(\mathcal{D}(G)\) \emph{perfect}?
  \item Under what conditions is \(\mathcal{D}(G)\) \emph{bipartite}?
  \item When does \(\mathcal{D}(G)\) possess an \emph{Eulerian} trail?
\end{itemize}

We present a complete characterization of those finite groups \(G\) (with \(|\pi(Z(G))|\ge2\)) for which \(\mathcal{D}(G)\) is bipartite. Additionally, we explore the perfectness and Eulerian properties of \(\mathcal{D}(G)\), shedding light on their structural and combinatorial implications.

\begin{theorem}\label{t:nilp}
Let $G$ be a finite nilpotent group with $\pi(G)\geq 3.$ Then $\mathcal{D}(G)$ is perfect if and only if $G\cong\mathbb{Z}_{p_1p_2p_3}.$    
\end{theorem}

To establish Theorem~\ref{t:nilp}, we begin with the following proposition:

\begin{proposition}\label{t:nilp_imperf}
Let \(G\) be a finite nilpotent group satisfying \(\pi(G) \ge 4\). Then the difference graph \(\mathcal{D}(G)\) is not perfect.
\end{proposition}

\begin{proof}
Let \(p_1, p_2, p_3, p_4\) be four distinct primes dividing \(\lvert G \rvert\). Since \(G\) is nilpotent, one can select elements
$a_1, \dots, a_5 \in G$
with the following orders:
\[
\mathrm{o}(a_1) = p_1 p_2,\quad
\mathrm{o}(a_2) = p_1 p_3,\quad
\mathrm{o}(a_3) = p_3 p_4,\quad
\mathrm{o}(a_4) = p_1 p_2 p_3,\quad
\mathrm{o}(a_5) = p_1 p_3 p_4.
\]
None of these \(a_i\) are isolated in \(\mathcal{D}(G)\), by Theorem~\ref{t:isol}.  
Notice that
$a_1 \sim a_2 \sim a_3 \sim a_4 \sim a_5 \sim a_1$
forms a 5-cycle an induced odd cycle of length 5 in \(\mathcal{D}(G)\). However, odd holes (induced odd cycles of length \(\ge 5\)) cannot appear in perfect graphs by the Strong Perfect Graph Theorem. Hence, \(\mathcal{D}(G)\) is not perfect.
\end{proof}

\begin{lemma}\label{t:nilp_imperf2}
Let \(G\) be a finite nilpotent group that is not cyclic and satisfies \(\pi(G)= 3\). Then the difference graph \(\mathcal{D}(G)\) is not perfect.
\end{lemma}

\begin{proof}
Let \(p_1, p_2, p_3\) be the distinct primes dividing \(\lvert G\rvert\). Since \(G\) is nilpotent and non-cyclic, there exists at least one prime, say \(p_1\), for which there are at least three distinct cyclic subgroups of order \(p_1\). Denote corresponding generators by \(a_1, a_2, a_3\) with \(\langle a_i\rangle \neq \langle a_j\rangle\) for \(i\neq j\). Pick elements \(b\) and \(c\) in \(G\) with \(\operatorname{o}(b) = p_2\) and \(\operatorname{o}(c) = p_3\), respectively. 

Then one can verify that the vertices
\[
a_1 b,\quad
a_2 b,\quad
a_1 b c,\quad
a_3 c,\quad
a_1 c
\]
form an induced cycle of length 5 in \(\mathcal{D}(G)\). Because induced odd cycles of length at least 5 (odd holes) are forbidden in perfect graphs by the Strong Perfect Graph Theorem, it follows that \(\mathcal{D}(G)\) cannot be perfect.
\end{proof}

\begin{lemma}\label{t:cyc_perf}
Let \(G\) be a cyclic group with \(\pi(G) = 3\). Then the difference graph \(\mathcal{D}(G)\) is perfect if and only if
$G \cong \mathbb{Z}_{p_1 p_2 p_3},$
where \(p_1, p_2, p_3\) are distinct primes.
\end{lemma}

\begin{proof}
First, assume \(G\) is cyclic and not of the form \(\mathbb{Z}_{p_1 p_2 p_3}\); then some prime square \(p_1^2\) divides \(\lvert G\rvert\). One constructs an induced 5-cycle in \(\mathcal{D}(G)\): select elements \(a_1, a_2, \dots, a_5\) such that
\[
\mathrm{o}(a_1) = p_1^2,\quad
\mathrm{o}(a_2) = p_1 p_2,\quad
\mathrm{o}(a_3) = p_2 p_3,\quad
\mathrm{o}(a_4) = p_1^2 p_2,\quad
\mathrm{o}(a_5) = p_1 p_2 p_3.
\]
They form an induced odd cycle of length 5; by the Strong Perfect Graph Theorem, such an odd hole prevents \(\mathcal{D}(G)\) from being perfect.
Now let \(G \cong \mathbb{Z}_{p_1 p_2 p_3}\), where \(p_1, p_2, p_3\) are distinct prime numbers. In this case, by Theorem~\ref{t:isol}, every non‑isolated vertex in the difference graph \(\mathcal{D}(G)\) has order of the form \(p_i p_j\) with \(i \ne j\).
Suppose, for contradiction, that \(\mathcal{D}(G)\) contains an induced odd cycle of length \(2k + 1\) with \(k \ge 2\):
\[
a_1 \sim a_2 \sim \cdots \sim a_{2k+1} \sim a_1.
\]
Without loss of generality, let \(\mathrm{o}(a_1) = p_1 p_3\). Then by Lemma~\ref{l:gens_null}, \(\mathrm{o}(a_2)\) must be either \(p_2 p_3\) or \(p_1 p_2\). 
If \(\mathrm{o}(a_2) = p_2 p_3\), then \(\mathrm{o}(a_3)\) is either \(p_1 p_3\) or \(p_1 p_2\).  
In the former case, \(\mathrm{o}(a_4)\) must be one of \(p_1 p_2\) or \(p_2 p_3\), which makes \(a_4\) adjacent to \(a_1\) (by Lemma~\ref{l:blowup}), contradicting the assumption that the cycle is induced.
Thus, no induced odd cycle of length $\geq5$ can exist in \(\mathcal{D}(G)\).
Moreover, the complement graph \(\overline{\mathcal{D}(G)}\) breaks into three disjoint cliques again precluding any induced odd cycle of length $\geq 5.$

\end{proof}

\begin{proof}[Proof of Theorem \ref{t:nilp}]
The proof follows from Proposition \ref{t:nilp_imperf} and Lemmas \ref{t:nilp_imperf2}, and \ref{t:cyc_perf}.    
\end{proof}

\subsection{Bipartite graph} 
\begin{theorem}\label{t:bip}
Let \(G\) be a finite group with \(\pi(Z(G))\ge 2\). Then the difference graph \(\mathcal{D}(G)\) is bipartite if and only if 
\(G \cong \mathbb{Z}_{p_i^{\alpha_i} p_j}\), where \(p_i, p_j\) are distinct primes and \(\alpha_i \ge 1\).
  %\item \(G \cong \mathbb{Z}_{p_i^{\alpha_i}}\).

\end{theorem}

\begin{proof}
Let \(G\) be a finite group with \(\pi(Z(G)) \ge 2\). We consider several cases:

\noindent\textbf{Case 1:} \(G\) is neither cyclic nor isomorphic to \(\mathbb{Z}_m \times Q_{2^n}\) (with \(\gcd(m, 2)=1\)).  
By Lemma~\ref{t:k3}, \(\mathcal{D}(G)\) contains a 3-cycle. Hence, it is \emph{not} bipartite (since bipartite graphs cannot contain odd cycles).

\noindent\textbf{Case 2:} \(G\) is cyclic with \(\lvert G \rvert = p_1^{\alpha_1} \cdots p_r^{\alpha_r}\) and \(r \ge 3\).  
Choose three distinct primes \(p_i, p_j, p_k\) dividing \(\lvert G \rvert\). There exist elements \(g_1, g_2, g_3 \in G\) with
\[
\mathrm{o}(g_1) = p_i p_j,\quad
\mathrm{o}(g_2) = p_j p_k,\quad
\mathrm{o}(g_3) = p_i p_k,
\]
which form a triangle in \(\mathcal{D}(G)\). Again, this odd cycle implies that \(\mathcal{D}(G)\) is not bipartite.

\noindent\textbf{Case 3:} \(\lvert G \rvert = p_1^{\alpha_1} p_2^{\alpha_2}\) with \(\alpha_1, \alpha_2 \ge 2\).  
Select five elements \(a_1, \dots, a_5 \in G\) with orders:
\[
\mathrm{o}(a_1) = p_1^2,\quad
\mathrm{o}(a_2) = p_1 p_2^2,\quad
\mathrm{o}(a_3) = p_1^2 p_2,\quad
\mathrm{o}(a_4) = p_2^2,\quad
\mathrm{o}(a_5) = p_1 p_2.
\]
These form a $5$-cycle in \(\mathcal{D}(G)\), proving that it is not bipartite.

\noindent\textbf{Case 4:} One exponent \(\alpha_i = 1\). Without loss of generality, let \(\alpha_2 = 1\), so
\[
\lvert G \rvert = p_1^{\alpha_1} p_2, \quad \alpha_1 \ge 1.
\]
If \(\alpha_1 = 1\), then by Theorem~\ref{t:isol}, \(\mathcal{D}(G)\) is an empty graph, hence bipartite.
If \(\alpha_1 > 1\), partition
\[
V(\mathcal{D}(G)) = V_1 \cup V_2,
\]
where
\[
V_1 = \{g \in G : \mathrm{o}(g) = p_1^t\},\qquad
V_2 = \{g \in G : \mathrm{o}(g) = p_1^t p_2\},
\]
for \(0 \le t \le \alpha_1\). There are no edges within \(V_1\) nor within \(V_2\), making \(\mathcal{D}(G)\) bipartite. 
\end{proof}

%%%%%%%%%%%%%%%%%%%%%%%%%%%%%%%%%%%%%%%%%%%%%%%%%%%%%%%%%%%%% Eulerian %%%%%%%%%%%%%%%%%%%%%%%%%%%%%%%%%%%%%%%%%%%%%%%%%%%%%%%%%%%%%%%%%%%

\subsection{Eulerian Properties of Difference Graphs}

In this section, we characterize all finite groups \(G\) for which the difference graph \(\mathcal{D}(G)\) is an even graph.

\begin{theorem}\label{t:euler}
For any finite group \(G\), the difference graph \(\mathcal{D}(G)\) is an even graph.
\end{theorem}

\begin{proof}
To show that \(\mathcal{D}(G)\) is an even graph, it suffices by Euler’s theorem to prove that each vertex has even degree. Since \(\mathcal{D}(G)\) is derived from the group structure, its connectedness is naturally inherited from adjacency definitions; so it remains only to verify that every vertex has even degree.

Let \(a\) be any vertex in \(\mathcal{D}(G)\).  Let $x_1, \dots, x_m\in  V(\mathcal{D}(G))\setminus\mathcal{G}_a$ be such that $a\sim x_i$ and $\mathcal{G}_{x_i}\cap \mathcal{G}_{x_j}=\emptyset,$ for all $i\neq j.$
By Lemma~\ref{l:blowup}, each neighbor \(x_i\) brings along all vertices in \(\mathcal{G}_{x_i}\) as neighbors of \(a\). Therefore, degree of \(a\) is
\[
\deg(a) = \sum_{i=1}^m \phi\bigl(\mathrm{o}(x_i)\bigr),
\]
where \(\phi\) is Euler’s totient function.

By Proposition~\ref{p-elts}, any non-isolated neighbor \(x_i\) must have \(\mathrm{o}(x_i) > 3\), so \(\phi\bigl(\mathrm{o}(x_i)\bigr)\) is even for each \(i\). Hence \(\deg(a)\), being a sum of even integers, is itself even.
Consequently, every vertex of \(\mathcal{D}(G)\) has even degree. This completes the proof that \(\mathcal{D}(G)\) is even graph.
\end{proof}

\section{Simple groups}
\label{sec:simple-grp}

Recall the class $\mathcal{C}$ of groups in which the difference graph is empty.
In this section we determine the non-abelian simple groups $G$ for which $\mathcal{D}(G)$ is empty. Since these groups
have trivial center, it follows from Theorem~\ref{t:induct} that to decide
whether $G\in\mathcal{C}$ it suffices to consider whether the proper subgroups
of such a group belong to $\mathcal{C}$. The only such groups are
$G=\mathrm{PSL}(2,q)$ for certain prime powers $q$, though it seems beyond the
reach of number theory to determine precisely which prime powers occur, or
even if there are finitely or infinitely many. 

We deal first with the groups $\mathrm{PSL}(2,q)$; this requires a lemma.

\begin{lemma}
The difference graphs of the cyclic group $\mathbb{Z}_n$ and the dihedral group
$D_{2n}$ are equal.
\end{lemma}

\begin{proof}
The elements of $D_{2n}\setminus\mathbb{Z}_n$ all have order~$2$, so by
Proposition~\ref{p-elts}(a), they are isolated, and
hence deleted when we construct the difference graph.
\end{proof}

\begin{proposition}\label{Prop:null-diff-simple-grp}
The simple group $\mathrm{PSL}(2,q)$ has difference graph empty if and only if
the following holds: setting $d=\gcd(q-1,2)$, each of $(q-1)/d$ and $(q+1)/d$ is
either a prime power or the product of two primes.
\end{proposition}

\begin{proof}
The list of maximal subgroups of $\mathrm{PSL}(2,q)$ was given by
Dickson~\cite[Chapter 22]{dickson} and is available in many places, for example Huppert~\cite[Chapter 2]{huppert}. They are among the
following (not all of these are maximal, or even occur in a given group):
\begin{enumerate}
\item Dihedral groups of order $2(q\pm1)/d$;
\item Semidirect products $V\rtimes\mathbb{Z}_{(q-1)/d}$;
\item Subgroups $\mathrm{PSL}(2,q_0)$ or $\mathrm{PGL}(2,q_0)$, where $q$ is a power of $q_0$;
\item $A_4$, $S_4$ or $A_5$.
\end{enumerate}
We consider these in turn.

(a) The difference graph for these is the same as that of the cyclic group
of order $(q-1)/d$; our earlier results show that it is empty if and only if
its order is either a prime power or the product of two primes.

(b) The affine group has trivial center, so we only have to consider the
elementary abelian group of order $q$ and the cyclic group of order 
$(q-1)/d$.

(c) These groups will not concern us since $q_0\pm1$ divides $q\pm1$.

(d) These three groups are easily seen to have empty difference graph. (The
groups $A_4$ and $A_5$ have all elements of prime order, while the maximal
subgroups of $S_4$ are $D_6$, $D_8$ and $A_4$.)
\end{proof}

\begin{remark}
The values of $q$ for which the difference graph of $\mathrm{PSL}(2,q)$ is empty are
the same as those for which the power graph of $\mathrm{PSL}(2,q)$ is a cograph
\cite[Proposition 8.7]{gog}. However, number theory does not yet provide the
tools to decide exactly which such groups occur.
\end{remark}

To deal with the remaining groups, we use Proposition~\ref{p:nulld_implies_cograph}, which shows that their power graphs must be cographs, together with
\cite[Theorem 1.3]{cmm}, which determines the non-abelian simple groups whose
power graphs are cographs: they are certain groups
$\mathrm{PSL}(2,q)$ and $\mathrm{Sz}(q)$ together with $\mathrm{PSL}(3,4)$.
So we only need to show that the difference graphs of Suzuki groups and 
$\mathrm{PSL}(3,4)$ are not empty, which we do by showing that these groups
contain the quaternion group $Q_8$ as a subgroup.

Let $q$ be an odd power of $2$, with $q>2$. The Suzuki group $\mathrm{Sz}(q)$
has Sylow $2$-subgroup $P$ of order $q^2$; its center is elementary abelian
of order $q$, and the remaining elements all have order $4$. Moreover, the
group has an automorphism of order $q-1$ permuting the involutions
transitively~\cite{higman:sz}. So each of the $q-1$ involutions is the square
of $q$ elements of order $4$; in particular, there exist $x$ and $y$ with
$x^2=y^2=z$ (say) and $z^2=1$. Thus $\langle x,y\rangle$ is a quaternion group.

Finally, consider $\mathrm{PSL}(3,4)$. We will only consider elements of
$2$-power order, so we may work in the matrix group $\mathrm{SL}(3,4)$.
Let
\[X_a=\begin{pmatrix}1&a&0\\0&1&a^{-1}\\0&0&1\\\end{pmatrix}\]
for each non-zero element $a$ in the field of $4$ elements. It is easily
checked that $\det(X_a)=1$ and that
\[X_a^2=Z=\begin{pmatrix}1&0&1\\0&1&0\\0&0&1\\\end{pmatrix}.\]
Again we have two such matrices $X_1$ and $X_\omega$ (where $\omega$ is a
cube root of $1$) such that $X_1^2=X_\omega^2=Z$, and so a quaternion group.
So we have proved:

\begin{theorem}\label{t:null-diff-graph-simple-grp}
The only non-abelian finite simple groups whose difference graphs are empty
are those in Proposition~\ref{Prop:null-diff-simple-grp}.
\end{theorem}

\section{Groups with a unique involution}
\label{sec:grp-with-unique-involution}

The class of groups with a unique involution has many occurrences in different
areas; for example, it contains the binary polyhedral groups
\cite[p. 82]{coxeter}, which feature in the celebrated observation by John
McKay connecting the regular polyhedra with the exceptional root systems
(discussed in~\cite[pp. 109f]{ade}).

Let $G$ be a group containing a unique involution $z$. Then the set $F$ of
elements of even order in $G$ contains at least half the elements of $G$
(since, of each pair $(g,gz)$, at least one belongs to $F$) and carries a
complete graph in the intersection power graph of $G$. Moreover, if $g\notin F$,
then $g$ is joined to at least the element $gz$ of $F$ in the power graph;
and, if the order of $G$ has a proper divisor $d$, then $g^dz$ is an element of
$F$ joined to $G$ in the difference graph of $G$. The difference graphs of at
least some of these groups, for example the binary polyhedral groups, should
repay investigation.

The structure of groups in $G$ is known \cite{bc}. We outline the argument.
First, note that the Sylow $2$-subgroups of $G$ also have a unnique involution,
and so are cyclic or generalized quaternion. Also, $Z=\langle z\rangle$ is a
normal (even central) subgroup of $G$, and $H=G/Z$ has cyclic or dihedral Sylow
$2$-subgroups. The structure of such groups is determined by Burnside's
transfer theorem~\cite[p. 155]{burnside} and the Gorenstein--Walter
theorem~\cite{gw}. If $O(H)$ is the largest odd-order normal subgroup of $H$,
then $H/O(H)$ is isomorphic to a subgroup of $\mathrm{PSL}(2,q)$
coontaining $\mathrm{PSL}(2,q)$, or $A_7$, or a Sylow $2$-subgroup of $H$.
Now a cohomological argument due to Glauberman, described in~\cite{bc}, shows
that, for any group $H$ with cyclic or dihedral Sylow $2$-subgroups, there is
a group $G$, unique up to isomorphism, which contains a unique involution $z$ 
and has the property that $G/\langle z\rangle$ is isomorphic to $H$.

\section{An example}
\label{sec:example}

In the paper \cite{bcdd}, the authors took various small finite simple groups
(with orders less than $10^4$), constructed the difference of the power graph
and the enhanced power graph, performed complete twin reduction on this graph,
and inspected the result. For most of the graphs tested, the results were not
too interesting; either the graph was reduced to a single vertex, or it had
a number of isomorphic connected components. For the Mathieu group $M_{11}$,
of order $7920$, an interesting graph was found; it was bipartite semiregular
on $385$ vertices with degrees $3$ and $4$ having diameter and girth $10$ and
automorphism group $M_{11}$.

We have performed some similar computations using the difference graph
$\mathcal{D}(G)$ defined in this paper. The results seem broadly similar.
In particular, for the Mathieu group $M_{11}$, we obtained a graph with the vertex set $V$ containing
$825$ vertices. It contains an independent set $A$ of size $165$, each of whose
vertices has four neigbours in the induced subgraph on the complement $B$, where $B=V\setminus A$.
This graph on $660$ vertices is regular with valency $5$; it has diameter~$9$,
girth~$3$, and automorphism group $M_{11}$. 
Further investigations would be interesting!

%\subsection*{Data Availability Statement} Availability of data and materials are not applicable.
%\subsection*{Conflict of interest}
%The authors declare that they have no conflict of interest.
%\bibliographystyle{amsplain}
%\bibliography{gen-inv-lcp.bib}

\section{Conclusion and open issues}
\label{sec:Conclusion-open-prob}

In this paper, we studied the difference graph $\mathcal{D}(G)$ of a finite group $G$. 
The investigation was primarily focused on identifying the isolated vertices in $\mathcal{D}(G)$, 
in particular, determining which groups yield a difference graph with no edges, 
as well as examining the connectedness and perfectness of such graphs. 
Several questions arising from this work suggest directions for further research.

For a finite group $G$ with $\pi(Z(G))\geq 2$, it was shown that $\mathcal{D}(G)$ is connected 
and has diameter at most $6$. 
However, for groups $G$ with $\pi(Z(G))\leq 1,$  $\mathcal{D}(G)$ may or may not be connected. 
This naturally leads to the following open problems.

\begin{question}\label{Q111}
If $\pi(Z(G))\leq 1,$ then what can be said about the connectivity of the difference graph $\mathcal{D}(G)?$ Can 
$\mathrm{diam}(\mathcal{D}(G))$ be greater than $6$? 
\end{question}

\begin{question}\label{Q112}
Complete the classification of isolated vertices of the difference graph.
\end{question}

\begin{question}\label{Q113}
Complete the classification of finite groups whose difference graph is perfect.
\end{question}

\begin{question}\label{Q114}
Complete the classification of finite groups whose difference graph is bipartite.
\end{question}

\section*{Acknowledgements} 
The authors would like to express sincere gratitude to the referee for careful
reading of the manuscript and for valuable comments and suggestions, which
have significantly improved the presentation and quality of this work. SB gratefully acknowledges the financial support provided by the NBHM research project (Reference No.~02011/29/2025NBHM(RP)/R\&DII/11951). He sincerely thanks the National Board for Higher Mathematics (NBHM), India, for funding this research. SB also acknowledges the excellent research environment and facilities provided by Dhirubhai Ambani University, Gandhinagar, Gujarat. 

\subsection*{Declaration of competing interest}
The authors have no relevant financial or non-financial interests to disclose.

\subsection*{Data availability}
Data sharing is not applicable to this article as no datasets were generated or analyzed 
during the present study.

\end{document}